\documentclass[journal, letterpaper]{IEEEtran}
\usepackage{amsfonts,amsmath,graphicx,amssymb,algorithm,algorithmic,epsfig,multirow}
\newcommand{\bl}{\boldsymbol}
\newcommand{\ml}{\mathcal}

\newtheorem{lemma}{Lemma}
\newtheorem{theorem}{Theorem}

\title{An Information-Theoretic Approach to PMU Placement in Electric Power Systems}
\author{Qiao Li, \IEEEmembership{Student Member, IEEE}, Tao Cui, \IEEEmembership{Student Member, IEEE}, Yang Weng, \IEEEmembership{Student Member, IEEE}, Rohit~Negi, \IEEEmembership{Member, IEEE}, Franz Franchetti \IEEEmembership{Member, IEEE}, and Marija D. Ili\'c, \IEEEmembership{Fellow, IEEE} 

\thanks{All authors are affiliated with the Department of Electrical and Computer Engineering, Carnegie Mellon University, Pittsburgh, PA, 15213 USA (email: qiaoli@cmu.edu, \{tcui, yweng, negi, franzf, milic\}@ece.cmu.edu).}}

\begin{document}

\maketitle

\begin{abstract}
  This paper presents an information-theoretic approach to address the phasor measurement unit (PMU)  placement problem in electric power systems. Different from the conventional `topological observability' based approaches, this paper advocates a much more refined, information-theoretic criterion, namely the mutual information (MI) between the PMU measurements and the power system states. The proposed MI criterion can not only include the full system observability as a special case, but also can rigorously model the remaining uncertainties in the power system states with PMU measurements, so as to generate highly informative PMU configurations. Further, the MI criterion can facilitate robust PMU placement by explicitly modeling probabilistic PMU outages. We propose a greedy PMU placement algorithm, and show that it achieves an approximation ratio of $(1-1/e)$ for any PMU placement budget. We further show that the performance is the best that one can achieve in practice, in the sense that it is NP-hard to achieve any approximation ratio beyond $(1-1/e)$. Such performance guarantee makes the greedy algorithm very attractive in the practical scenario of multi-stage installations for utilities with limited budgets. Finally, simulation results demonstrate near-optimal performance of the proposed PMU placement algorithm.
\end{abstract}


\section{Introduction}
\label{sec:intro}

\IEEEPARstart{S}{ynchronized} measurement technology (SMT) has been widely recognized as an enabler of the emerging real-time wide area monitoring, protection and control (WAMPAC) systems \cite{terizija11, bakken11}. Phasor measurement unit (PMU), being the most advanced and accurate instrument of SMT, plays a critical role in achieving key WAMPAC functionalities \cite{phadke10}. With better than one microsecond global positioning system (GPS) synchronization accuracy, the PMUs can provide highly synchronized, real-time, and direct measurements of voltage phasors at the installed buses, as well as current phasors of adjacent power branches. Such measurements are vital for the efficient and reliable operations of the power systems by both improving the Situational Awareness (SA) of the grid operators, and facilitating synchronized and just-in-time (JIT) automated control actions \cite{sun07, overbye10} .

Given the critical role of PMUs for the power system, it is important that these instruments are installed at carefully chosen buses, so as to maximize the `information gain' on the system states, and achieve desired functionalities efficiently. Currently, there is a significant performance gap between the existing research and the desired `informative' PMU configuration \cite{rice06, li11}. One particular reason is that most researches center around the \emph{topological observability} criterion, which only specifies that power system states should be uniquely estimated using minimum number of PMU measurements \cite{baldwin93}. Based on such criterion, many solutions were proposed, such as the ones based on mixed integer programming \cite{xu04, gou08}, binary search \cite{chakrabarti08}, and metaheuristics \cite{milosevic03, aminifar09}. While it is true that all PMU configurations can monitor the power system states with similar accuracy once the system becomes fully observable, these PMU placement approaches can yield quite suboptimal results for the important and current situation, where the number of installed PMUs is far from sufficient to achieve full system observability. The reason is as follows. Firstly, the `observability' criterion is very coarse, which specifies the information gained on the system states as \emph{binary}, i.e., either observable or non-observable. Such crude approximation essentially assumes that the states at different buses are completely independent (with exceptions for buses with zero injection), in that the knowledge of the state of a bus has zero information gain on the state of any other bus, as long as that bus is not `observable'. This is clearly not the case for power systems, where the system states exhibit high correlations, due to the fundamental physical laws, such as KVL and KCL. Secondly, the observability approaches neglect important parameters of the power system, such as transmission line admittances, by focusing only on the binary connectivity graph. In this sense, if zero injection are not considered, the current researches is essentially the classic `dominating set' problem \cite{vazirani04}, where a subset of buses in the system are selected, so that every bus is either in the subset, or a neighbor of the subset. Such over-simplification of the power system is very likely to result in suboptimal design and significant performance loss. For example, it has been shown in \cite{rice06, li11} that PMU configurations can have large influence on the accuracy of the state estimation, even though the observability result stays the same.

To overcome the performance limitation of current approaches, in this paper, we advocate a much more refined, information-theoretic criterion to generate highly informative PMU placement configurations. Specifically, we rigorously model the `information gain' achieved by the PMUs states as the Shannon \emph{mutual information} (MI) \cite{cover06} between the PMU measurements and the power system states. The MI criterion is very popular in the statistics and machine learning literature \cite{mackay92, krause08}, which has found many applications in sensor placement problems. For power systems, we will show that it can include the `topological observability' by current researches as a special case. Not only this, the MI criterion can also incorporate probabilistic PMU failures, to facilitate robust PMU placement configurations. A related work is the entropy based approach \cite{kamwa02}, where the PMU buses were selected, so that the `information content' of the PMU response signals from certain transient-stability program can be maximized. In \cite{kamwa02}, the `information' is represented by the norm of the entropy matrix of PMU response signals. Compared to their method, the MI criterion in this paper directly models the uncertainties in the power system states. Further, the MI criterion is rigorous, based on the analytical DC model of the power system.



As a second contribution of this paper, we present a greedy PMU placement algorithm, and show that it can achieve $(1-1/e)$ of the optimal information gain for any PMU budget $K$. We further prove that the approximation ratio is the best that one can achieve in practice, by showing that it is NP-hard to approximate with any factor larger than $(1-1/e)$. Compared to existing approaches, the greedy algorithm can not only achieve the best performance guarantee, but also can be easily extended to large-scale power systems. Further, the greedy PMU placement is very attractive in the practical scenario of multi-stage PMU installation, where utilities prefer to install the PMUs over a horizon of multiple periods, due to limited budgets \cite{nuqui05}. In such cases, utilities can simply adopt the greedy placement strategy, as the $(1-1/e)$ approximation ratio holds for any $K$. On the other hand, existing multi-stage methods \cite{nuqui05, aminifar11} may incur significant performance loss if the multi-period budget changes unexpectedly.

The remaining of this paper is organized as follows. Section \ref{sec:model} describes the power system and measurement models, and Section \ref{sec:pmu} formulates the optimal PMU placement problem. Section \ref{sec:greedy} proposes the greedy PMU placement algorithm and analyzes its performance, and Section \ref{sec:numerical} demonstrates the numerical results. Finally, Section \ref{sec:conclusion} concludes this paper.

\section{System Model}
\label{sec:model}

In this section, we formulate a Gaussian Markov random field (GMRF) model \cite{he11} for the system states, and describe the measurement models.

\subsection{GMRF Model for Phasor Angles}

A DC power flow model \cite{grainger94} is assumed in this paper, where the power injection $P^{\text{inj}}_i$ at bus $i$ can be expressed as follows:
\begin{equation}
  P^{\text{inj}}_i = B_{ii}\theta_i + \sum_{j\in\ml N_i} B_{ij}\theta_j
  \label{eqn:dc}
\end{equation}
In above, $\ml N_i$ is the set of neighboring buses of $i$, $B_{ij}$ is imaginary part of the nodal admittance matrix $Y$, and $\theta_i$ is the voltage phasor angle at bus $i$. The uncertainties in the power injection vector $\bl P^{\text{inj}}$ can often be approximated as Gaussian by existing stochastic power flow methods \cite{schellenberg05}. In this paper, we assume that $\bl P^{\text{inj}}$ is distributed as $\ml N(\bl \mu, \Sigma)$. Denote bus 0 as the slack bus. We are only interested in the states at non-reference buses, as the angle of the slack bus can be uniquely specified by the non-reference bus angles, due to the law of conservation of energy. Write the non-reference bus angles in vector form as $\bl \theta=(\theta_1, \theta_2, \ldots, \theta_N)^T$. Note that the system states $\bl \theta$ are highly correlated statistically, due to the DC model in (\ref{eqn:dc}). Formally, the dependency of these variables are described by the following theorem:

\begin{theorem}
  Assume the power system is fully connected. Under the DC model, $\bl \theta$ forms a GMRF with mean $B^{-1}{\bl \mu}$ and covariance matrix $B^{-1}{\Sigma}B^{-1}$.
  \label{thm:gmrf}
\end{theorem}
\begin{IEEEproof}
  Since the power system is fully connected, the matrix $B$ is invertible \cite{krumpholz80}. Thus, the states can be calculated as $\bl \theta = B^{-1}\bl P^{\text{inj}}$, from which the theorem follows.
 \end{IEEEproof}

\begin{figure}[!t]
  \centering    
  \epsfig{file=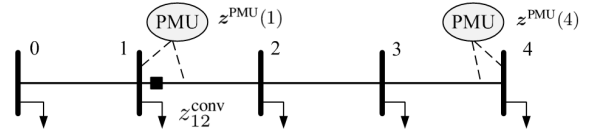,width=1\linewidth,clip=}\\
  \caption{The one line diagram of a power system with 5 buses. The square node represents the measurement of active power flow from bus 1 to 2. Two PMUs are installed at bus 1 and 4.}\label{fig:sample}
\end{figure}

Fig. \ref{fig:sample} illustrates a 5-bus power system, with its GRMF model shown as the shaded region in Fig. \ref{fig:graph}. In this case, the GRMF is formed by connecting two-hop neighbors of the buses in the original power system, as the power injections are assumed to be independent. We next describe the PMU measurement model.

\subsection{Measurement Model}

\subsubsection{Conventional Measurements}

As the DC model is assumed in this paper, the conventional measurements only include the real power injection $\{z^{\text{conv}}_i\}$ and real power flow $\{z^{\text{conv}}_{ij}\}$. Under the DC model, these measurements can be described as follows:
\begin{eqnarray}
  z_i^{\text{conv}} &=& B_{ii}\theta_i+\sum_{j\in\ml N_i} B_{ij} \theta_j+e^{\text{conv}}_i\\
  z^{\text{conv}}_{ij} &=& B_{ij}(\theta_j-\theta_i) + e^{\text{conv}}_{ij}
\end{eqnarray}
where $e^{\text{conv}}_i$ and $e^{\text{conv}}_{ij}$ are measurement noises, which are distributed as $\ml N(0, \kappa^{\text{conv}}_i)$ and $\ml N(0, \kappa^{\text{conv}}_{ij})$, respectively.

\subsubsection{PMU Measurements}

A PMU placed at bus $i$ can measure both the voltage at bus $i$ and the currents of selected incident branches. This implies that the phasor angles of corresponding adjacent buses can also be directly calculated. Thus, we assume the following equivalent PMU measurement model. We associate a PMU placed at bus $i$ with a vector $\bl z^{\text{PMU}}(i)$, such that
\begin{eqnarray}
  z_i^{\text{PMU}}(i) &=& \theta_i + e^{\text{PMU}}_i\label{eqn:zpmui}\\  
  z^{\text{PMU}}_{ij}(i) &=& (\theta_i-\theta_j) + e^{\text{PMU}}_{ij}, \forall j\in \ml P_i\label{eqn:zpmuij}
\end{eqnarray}
where $e^{\text{PMU}}_i$ and $e^{\text{PMU}}_{ij}$ are measurement noises with distribution $\ml N(0, \kappa^{\text{PMU}}_i)$ and $\ml N(0, \kappa^{\text{PMU}}_{ij})$, respectively. $\ml P_i\subseteq \ml N_i$ is a subset of neighbors of bus $i$. This is because of the PMU channel limits, which implies that only a subset of adjacent branches can be monitored. The variances $\kappa^{\text{PMU}}_{i}$ and $\kappa^{\text{PMU}}_{ij}$ depend on various sources of uncertainties, such as the GPS synchronization, instrument transformers, A/D converters and cable parameters, which can be estimated appropriately \cite{chakrabarti09}. Further, PMU failures can be modeled by assuming that each current measurement $z_{ij}^{\text{PMU}}(i)$ outputs a failure message with probability $1-a^{\text{PMU}}_{ij}$, and similarly, each voltage measurement fails with probability $1-a^{\text{PMU}}_i$, where $a^{\text{PMU}}_{ij}$ and $a^{\text{PMU}}_i$ are the availability of the current and voltage measurements, respectively.

\begin{figure}[!t]
  \setlength{\abovecaptionskip}{-1em}
  \centering    
  \epsfig{file=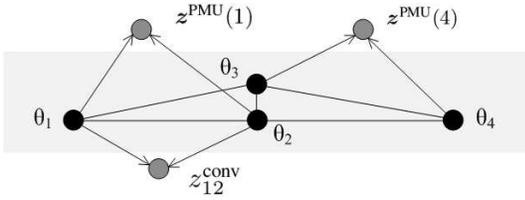,width=1\linewidth,clip=}
  \caption{The probabilistic graphical model for the non-reference bus angles in the power system in Fig. \ref{fig:sample}. The shaded region illustrates the GMRF for system states. The power injections are assumed to be independent.}
  \label{fig:graph}
\end{figure}

As an illustration, Fig. \ref{fig:graph} represents the probabilistic graph model for the measurement configuration in Fig. \ref{fig:sample}, where the gray nodes represents the measurement variables. From the figure, it is clear that the information gain of PMU measurements depends heavily on the placement buses. This will be formalized in the next section by the optimal PMU placement problem.

\section{Optimal PMU Placement}
\label{sec:pmu}

The placement configuration of PMUs should be highly `informative' to effectively monitor power system states. In this paper, we advocate an information-theoretic criterion to quantitatively assess the `information gain' that can be obtained from the PMU measurements. Specifically, we model the uncertainties in the system states as the Shannon entropy \cite{cover06}:
\begin{equation}
  H(\bl \theta)=\sum_{\bl \theta} -p(\bl \theta)\log p(\bl \theta)
\end{equation}
where $p(\bl \theta)$ is the probability mass function of $\bl \theta$. In this paper, we assume that the entropies of all variables are calculated after quantization with a sufficiently small step size $\Delta$. This is motivated by the fact that the phasor angles in power systems are observed by meters with finite accuracy. Denote $\ml S$ as the set of PMU configurations, where each element $s=\{i, \ml P_i\}$ in the set $\ml S$ corresponds to a candidate PMU configuration as in (\ref{eqn:zpmui}) and (\ref{eqn:zpmuij}). Note that our model is very general, which can be used to model PMU channel limits. The information gain of the PMU configuration $\ml S$ can be assessed by the entropy reduction due to PMU measurements $\bl z^{\text{PMU}}(\ml S)$:
\begin{equation}
  I(\bl \theta; \bl z^{\text{PMU}}(\ml S))=H(\bl \theta)-H(\bl \theta|\bl z^{\text{PMU}}(\ml S))
  \label{eqn:mi}
\end{equation}
where $I(\bl \theta; \bl z^{\text{PMU}}(\ml S))$ is the Shannon mutual information (MI) between the PMU measurements and the power system states, and $H(\bl \theta|\bl z^{\text{PMU}}(\ml S))$ is the conditional entropy. Finally, when conventional measurements are considered, the uncertainty reduction corresponds to the conditional MI:
\begin{equation*}
  I(\bl \theta; \bl z^{\text{PMU}}(\ml S)|\bl z^{\text{conv}})=H(\bl \theta|\bl z^{\text{conv}})-H(\bl \theta|\bl z^{\text{PMU}}(\ml S), \bl z^{\text{conv}})
\end{equation*}
The MI criterion is widely adopted in the machine learning literature to generate highly informative sensor placement configurations \cite{mackay92}. For the PMU placement problem in power systems, we claim that the MI criterion is also intimately related to the power system observability and state estimation accuracy, which we elaborate as follows:

{\it 1) (Observability):} The proposed MI criterion can include the observability criterion as a special case. To see this, assume no PMU failure and PMU measurement noise. We claim that the maximum information gain is achieved if and only if the power system is completely observable from PMU measurements. This can be clearly observed from (\ref{eqn:mi}), where the MI function is maximized if and only if $H(\bl \theta|\bl z^{\text{PMU}}(\ml S))=0$, in which case the system states are deterministic given the PMU measurements $\bl z^{\text{PMU}}(\ml S)$.

{\it 2) (State Estimation Error):} The proposed MI function is also intimately related to the minimization of the state estimation error. In fact, the conditional entropy $H(\bl \theta|\bl z^{\text{PMU}})$ can then be expressed as follows:
    \begin{equation}
      H(\bl \theta|\bl z^{\text{PMU}}) \approxeq \log\det \text{Cov}(\bl e)+C-N\log\Delta
      \label{eqn:mmse}
    \end{equation}
    where $C$ is a constant, and $\bl e=\bl \theta-\hat{\bl \theta}_{\text{MMSE}}(\bl z^{\text{PMU}})$ is the estimation error of the Minimum Mean Square Error (MMSE) estimation of $\bl \theta$ given $\bl z^{\text{PMU}}$, and $\text{Cov}(\cdot)$ is the covariance matrix. The error  in the above approximation is due to the quantization, and goes to zero as the quantization step size $\Delta\rightarrow 0$. Thus, it is clear that the maximization of MI can also lead to minimizing the MMSE error $\bl e$ of the power system states. Intuitively, $H(\bl \theta|\bl z^{\text{PMU}})$ specifies how `peaked' the distribution of estimated power system states behaves around the mean value. In the statistics literature, such criterion is referred to as the `D-optimality' \cite{mackay92}.

We are now ready to formulate the optimal PMU placement problem. Assume that there are a total of $K$ PMUs to be installed in the power system. The goal is to choose a subset of PMU configurations $\ml S^{\star}$ from a set of candidate PMU configurations, such that 
\begin{equation}
  \ml S^\star \in \arg\max_{|\ml S|\leq K} F_i({\ml S}), \qquad i = 1, 2.\
  \label{eqn:opp}
\end{equation}
The objective functions are illustrated as follows:

{\it 1) (PMU Measurements Only)} In this case, the objective function associated with a PMU placement set $\ml S$ is
\begin{eqnarray}
  F_1(\ml S) &=& {1\over T}\sum_{t=1}^T I_t(\bl \theta; \bl z^{\text{PMU}}({\ml S}))
  \label{eqn:F1}
\end{eqnarray}
where the dependence on time index $t$ is because the power system states $\bl \theta$ have time-dependent distribution $\ml N(B^{-1}\bl \mu_t, B^{-1}\Sigma_t B^{-1})$, due to the changes in real power injections over time. Thus, the objective function in (\ref{eqn:F1}) describes the `time averaged information gain' about the power system state over a time period of interest (such as one day).

{\it 2) (With Conventional Measurements)} When conventional measurements are considered, the objective function should be replaced with the conditional MI function, as follows:
\begin{eqnarray}
  F_2(\ml S) &=& {1\over T}\sum_{t=1}^T I_t(\bl \theta; \bl z^{\text{PMU}}({\ml S})| \bl z^{\text{conv}})
  \label{eqn:F2}
\end{eqnarray}
Note that it is possible that the time scales can be different in both cases, as the conventional measurements can have much slower sampling rate (on the order of minutes) than PMU measurements. Having formulated the optimal PMU placement problem, we will discuss the solutions in the next section.

\section{Greedy PMU Placement}
\label{sec:greedy}

It is highly desired that the PMUs are optimally placed in the power system. However, the optimal solution is very hard to obtain, as the optimal PMU placement problem is NP-complete \cite{brueni05}. In this section, we propose a greedy PMU placement algorithm, and show that it can achieve the optimal performance guarantee among the class of polynomial time algorithms.

\subsection{Hardness Result}

Before presenting the greedy algorithm, we first demonstrate the hardness result. We extend the hardness result in \cite{brueni05}, by showing that the optimal PMU placement problem is not only NP-hard to solve, but also NP-hard to approximate beyond the approximation ratio of $(1-1/e)$:
\begin{theorem}
  Unless $P=NP$, there is no polynomial time algorithm for the optimal PMU placement problem in (\ref{eqn:opp}) with better approximation ratio than $(1-1/e)$.
  \label{thm:hardness}
\end{theorem}
\begin{IEEEproof}
  See in Appendix \ref{apdx:hardness}.
\end{IEEEproof}
We next propose a greedy PMU placement algorithm, which can achieve the $(1-1/e)$ approximation ratio.

\subsection{Greedy PMU Placement}

The greedy PMU placement algorithm is shown in Algorithm \ref{alg:greedy}. Compared to the optimal placement, the greedy algorithm has low complexity, and is easy to implement in large-scale systems. In each step, the algorithm chooses the next candidate PMU configuration that can achieve the largest `marginal information gain', where the objective function $F(\cdot)$ can be chosen as either $F_1$ or $F_2$, depending on whether conventional measurements are included.
\begin{algorithm}[H]
  \caption{{Greedy PMU Placement}} \label{alg:greedy}
  \begin{algorithmic}[1]
    \STATE {\bf Initialize: } $\ml S\gets\emptyset$;
    \FOR{$k=1$ to $K$}
    \STATE $\ml S \gets \ml S\cup\{s^\star\}$, where $s^{\star}$ solves the following:
    \begin{equation}
      s^{\star} = \arg\max_{s\not\in \ml S} F(\ml S\cup\{s\});
      \label{eqn:update}
    \end{equation}
    \ENDFOR
    \RETURN{$\ml S$}
  \end{algorithmic}
\end{algorithm}
The next theorem shows that the greedy algorithm can achieve the largest approximation ratio of $(1-1/e)$.
\begin{theorem}
  The greedy PMU placement in (\ref{alg:greedy}) can achieve an approximation ratio of $(1-1/e)$ for both objective functions $F_1(\cdot)$ and $F_2(\cdot)$.
  \label{thm:greedy}
\end{theorem}
\begin{IEEEproof}
  The proof is obtained by identifying a key property, \emph{submodularity}, of the PMU placement problem. Detailed proof is in Appendix \ref{apdx:greedy}.
\end{IEEEproof}

We have the following remarks:

\subsubsection{Optimality} 

Based on Theorem \ref{thm:hardness} and \ref{thm:greedy}, we claim that the greedy algorithm can achieve the \emph{best performance guarantee that is possible}. Further, compared to methods such as mixed integer programming \cite{xu04, gou08}, binary search \cite{chakrabarti08}, or metaheuristics \cite{milosevic03, aminifar09}, the greedy algorithm is not only the best in performance guarantees, but also can be easily implemented in large-scale systems, due to the low computation complexity.

\subsubsection{Multi-stage Installation}

The greedy algorithm is very attractive in the case of multi-stage installations, where the utilities plan to install the PMUs over a horizon of multiple years, due to the limited (and possibly uncertain) budgets. In such scenarios, the greedy algorithm can always achieve an approximation ratio of $(1-1/e)$ for any given $K$, whereas fixed multi-stage planning algorithms may suffer from substantial performance loss when the budget $K$ changes unexpectedly.

Having formulated the greedy PMU placement algorithm and proved its optimality results, we will test it against other methods in standard IEEE test systems in the next section.

\section{Numerical Results}
\label{sec:numerical}

\begin{table}[!t]
\renewcommand{\arraystretch}{1.3}
\caption{Computation Time}
\label{table:computation} 
\centering
\begin{tabular}{c|c}
\hline
\bfseries Test Systems & \bfseries Time (Seconds)\\
\hline\hline
\begin{minipage}[t]{3.6cm}
  \centering
  IEEE 14-bus (PMU Only)
\end{minipage} & 0.7093  \\
\hline
\begin{minipage}[t]{3.6cm}
  \centering
  IEEE 14-bus (Conventional)
\end{minipage} &   2.5543    \\
\hline
\begin{minipage}[t]{3.6cm}
  \centering
  IEEE 57-bus (PMU Only)
\end{minipage} &        $2.9141\times 10^3$       \\
\hline
\begin{minipage}[t]{3.6cm}
  \centering
  IEEE 57-bus (Conventional)
\end{minipage} &         $3.4556\times 10^4$      \\
\hline
\end{tabular}
\end{table}

This section demonstrates the performance of the proposed greedy PMU placement algorithm, and compare it with the PMU placement results in the literature.

\subsection{System Description}

\begin{figure*}[!t]
  \centering
  \begin{tabular}{c c}
  \psfig{figure=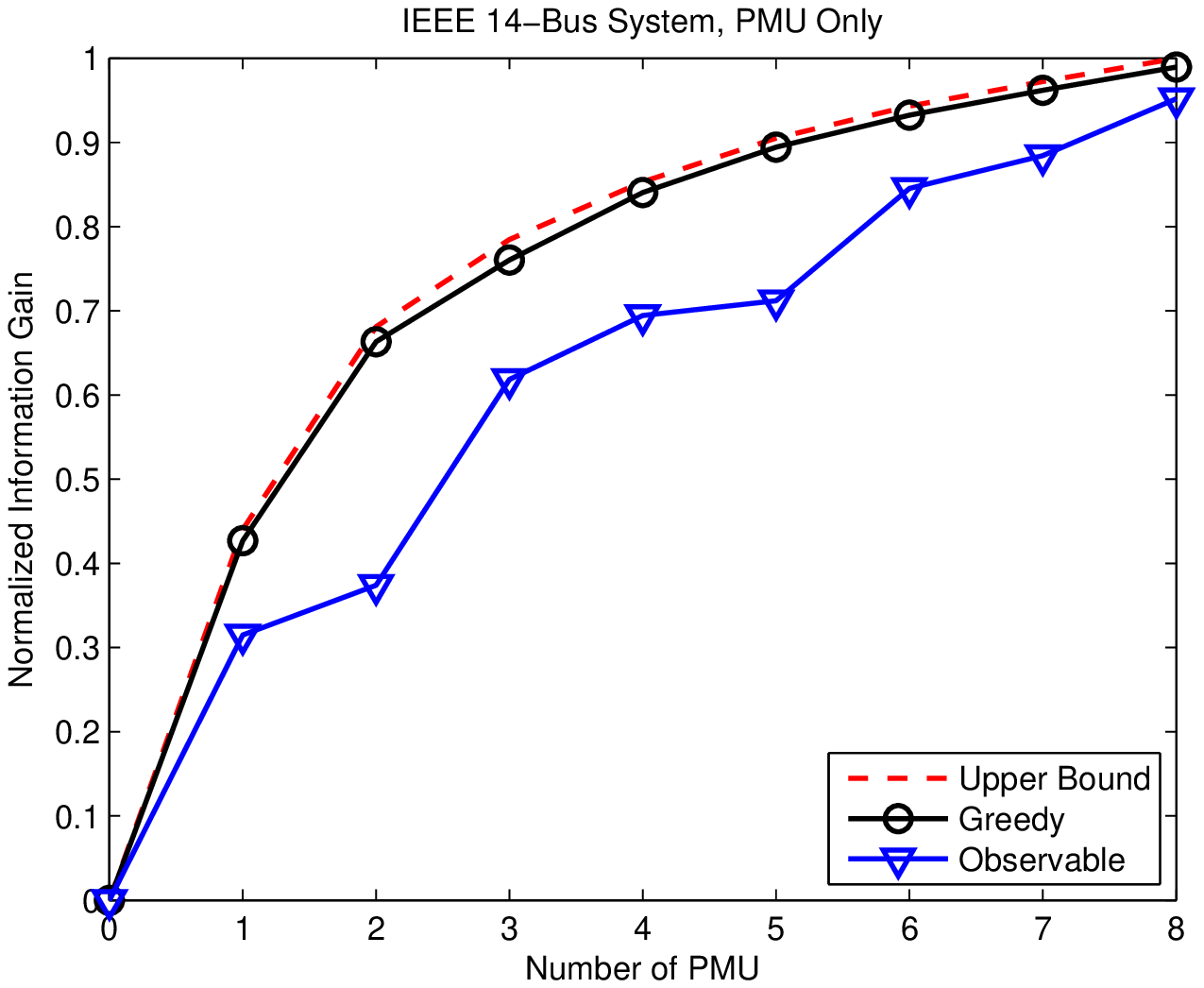,width=3.1in} & \psfig{figure=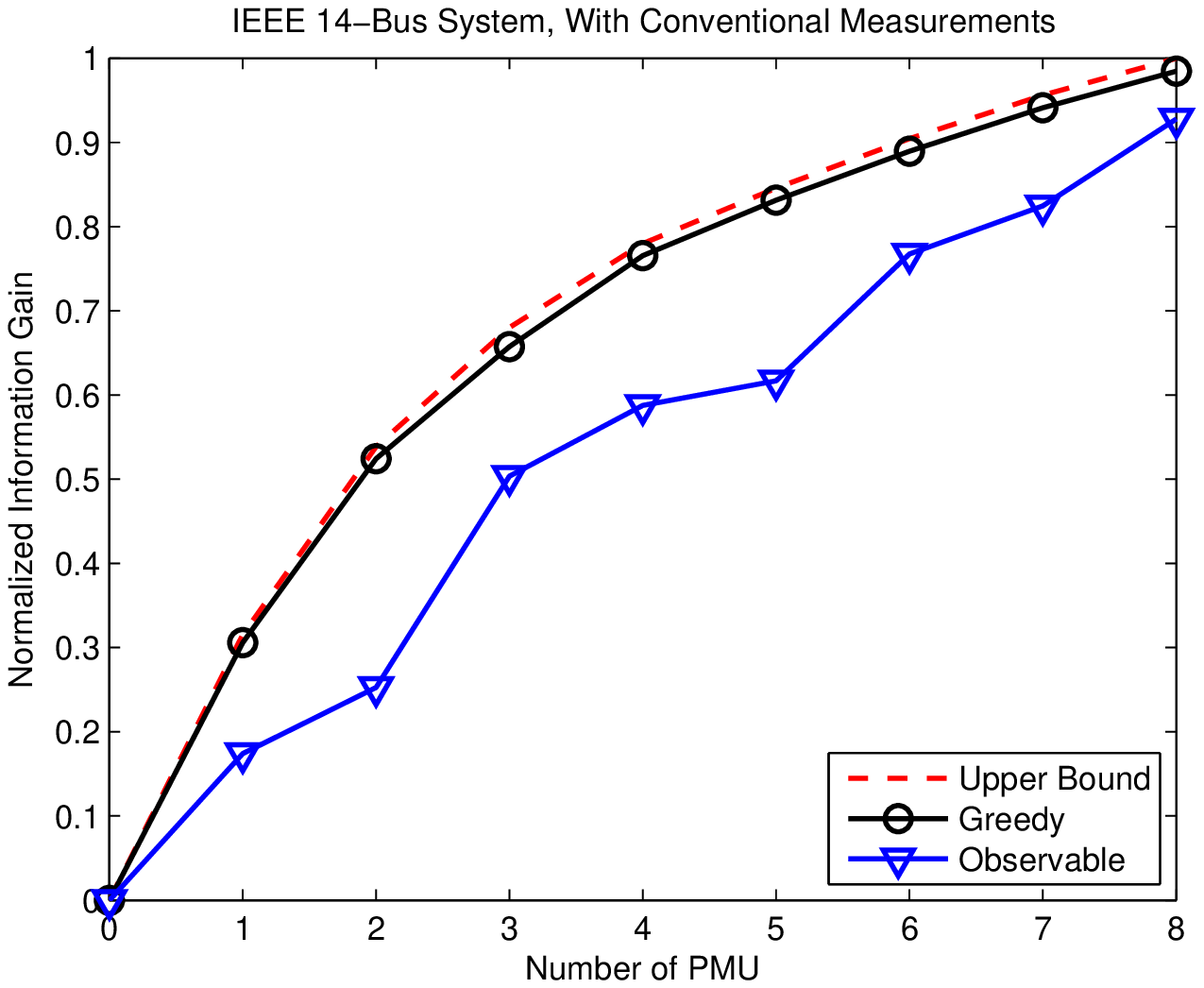,width=3.1in}\\
  (a) & (b)
  \end{tabular}
  \caption{Normalized information gain of different PMU placement schemes in the IEEE 14-bus system. (a) with PMU measurements only and (b) with conventional measurements.}\label{fig:ieee14}
\end{figure*}

In the simulation, real power injections are assumed to be normally distributed, and independent across different buses \cite{schellenberg05}. For each bus, the standard deviation of the real power injection is assumed to be $10\%$ of the mean value. For each time slot, the mean real power injections at all buses are obtained by properly scaling the case profile description of the standard test system \cite{christie99}. Thus, the MI function at different time slots are only different by a multiplicative factor. Each PMU measurement is assumed to fail independently with probability $0.03$ \cite{aminifar11}. The standard error of all PMU measurements are assumed to be $0.02^\circ$, whereas the standard error of conventional measurements are assumed to be $0.57^\circ$. Finally, for comparison purposes, we consider the `topological observability' based PMU placement configurations in \cite{xu05}, which are obtained based on mixed integer programming. All simulation results are obtained with MATLAB on an Intel Xeon E2540 CPU with 8GB RAM. The computation time is shown in Table \ref{table:computation}.

\begin{table}[!t]
\renewcommand{\arraystretch}{1.3}
\caption{PMU Locations for IEEE 14-Bus System}
\label{table:ieee14}
\centering
\begin{tabular}{c|c|c|c|c}
\hline
\bfseries $K$ & \bfseries 
\begin{minipage}[t]{1.7cm}
  \centering
 Optimal
 (PMU Only)
\end{minipage} & \bfseries 
\begin{minipage}[t]{1.7cm}
  \centering
 Greedy \\
 (PMU Only)
\end{minipage} & \bfseries 
\begin{minipage}[t]{1.6cm}
  \centering
 Optimal \\
 (Conventional)
\end{minipage} & \bfseries 
\begin{minipage}[t]{1.6cm}
  \centering
 Greedy\\
 (Conventional)
\end{minipage}\\
\hline\hline
1 & 4 & 4 & 6 & 6    \\
\hline
2 & 4, 13 & 4, 13 & 4, 13 & 6, 9       \\
\hline
3 & 4, 6, 9 & 4, 9, 13 & 4, 6, 14  & 4, 6, 9  \\
\hline
4 & 4, 6, 9, 13 & 4, 6, 9, 13 & 4, 6, 9 ,13 & 4, 6, 9, 13\\
\hline
\end{tabular}
\end{table}

\subsection{IEEE 14-Bus System}

\subsubsection{PMU Measurements Only}

\begin{figure}[!t]
  \centering
  \includegraphics[width=3.2in]{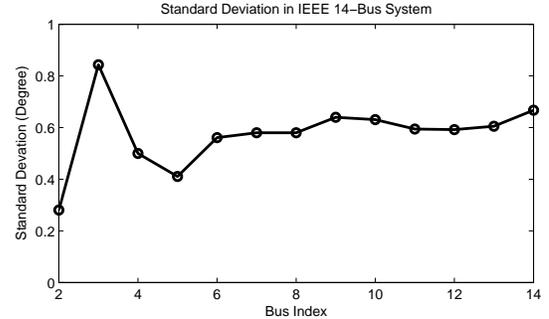}
  \caption{Standard deviation of voltage angles in the IEEE 14-bus system.}\label{fig:sd14}
\end{figure}

For this case, the optimal PMU locations are calculated by an exhaustive search among all possible configurations. The PMU locations for both optimal placement and greedy placement for $K\leq 4$ are shown in Table \ref{table:ieee14}. From the table, one can observe that the optimal PMU configuration can change significantly over different placement budget $K$. For example, when $K=2$, the optimal placement is $\{4, 13\}$, whereas for $K=3$, the optimal placement is $\{4, 6, 9\}$. On the other hand, the greedy placement configuration $\ml S_g(K)$ always satisfy $\ml S_g(K-1)\subset \ml S_g(K)$. Thus, the greedy placement strategy is robust against the uncertainties in the placement budget $K$. To verify the placement results, we plot the standard deviations of the phasor angles at all non-reference buses in Fig. \ref{fig:sd14}. From the figure, one can observe that the state at bus $3$ has the largest variance. However, bus $3$ is not chosen as the first PMU bus, since, intuitively, it is connected to only two neighbors in the system (see the topology in \cite{christie99}). Instead, bus $4$ is chosen, since it has five neighbors. 


Fig. \ref{fig:ieee14} (a) shows the normalized information gain for the IEEE 14-bus system with only PMU measurements. In the figure, the `Upper Bound' curve is computed by the optimal PMU configuration assuming no PMU failure. Thus, it overestimates the information gain on the system states. One can easily observe the near optimal performance of the greedy PMU placement strategy, in that the `Greedy' curve is very close to the `Upper Bound'. Further, the greedy algorithm has a significant improvement on information gain compared to the conventional `observability' based approach. For example, for $K=3$, the improvement is around $20\%$ as compared to the observability based placement \cite{xu05}. This clearly demonstrates the performance loss associated with the coarse observability based criterion. Finally, one can observe from the `Upper Bound' curve that the maximum information gain has a `diminishing marginal return' property, in that the marginal information gain tends to decrease as the number of installed PMUs in the power system grows. This also confirms the submodularity of the MI objective function.

%

\begin{figure*}[!t]
  \centering
  \begin{tabular}{c c}
  \psfig{figure=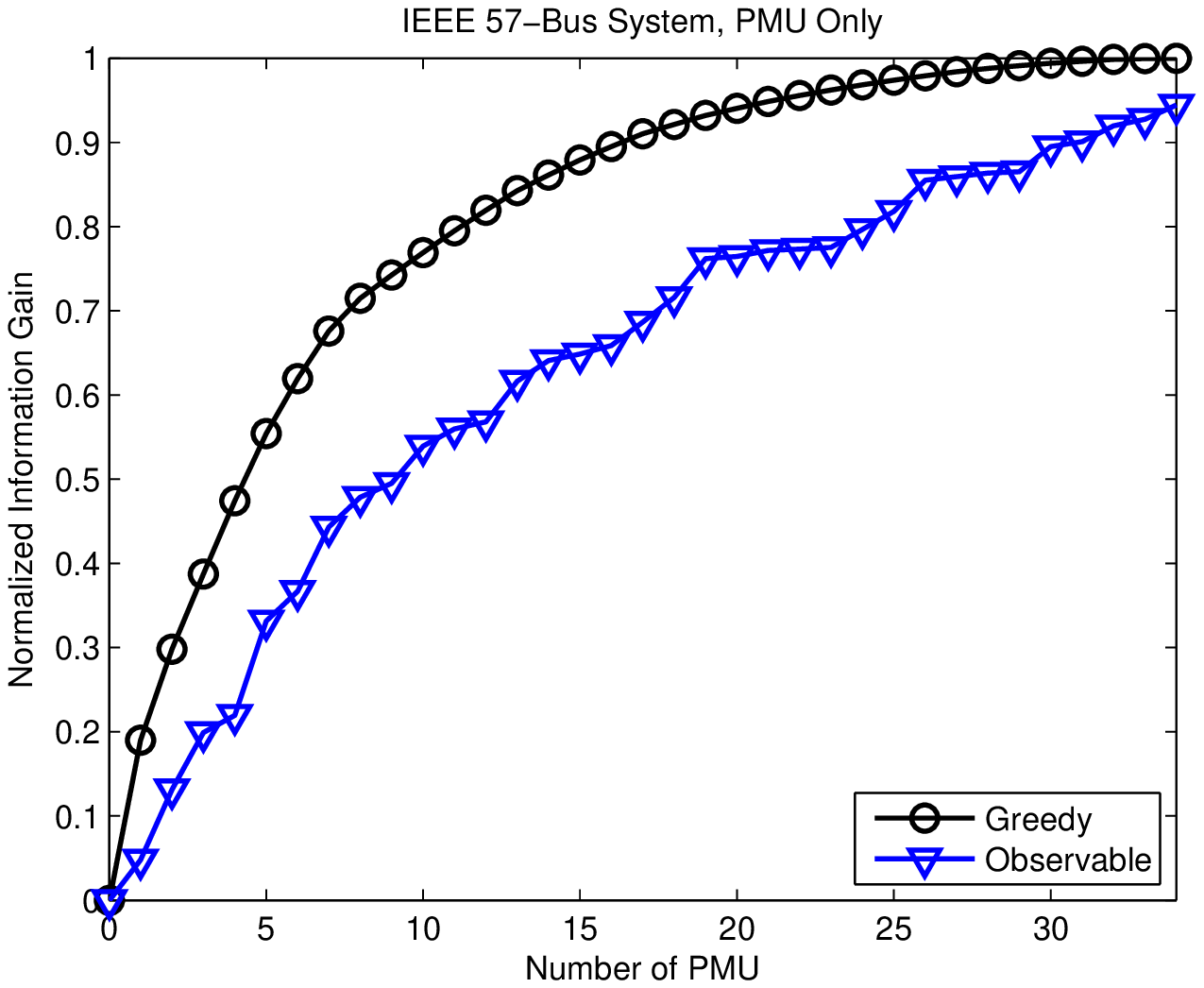,width=3.1in} & \psfig{figure=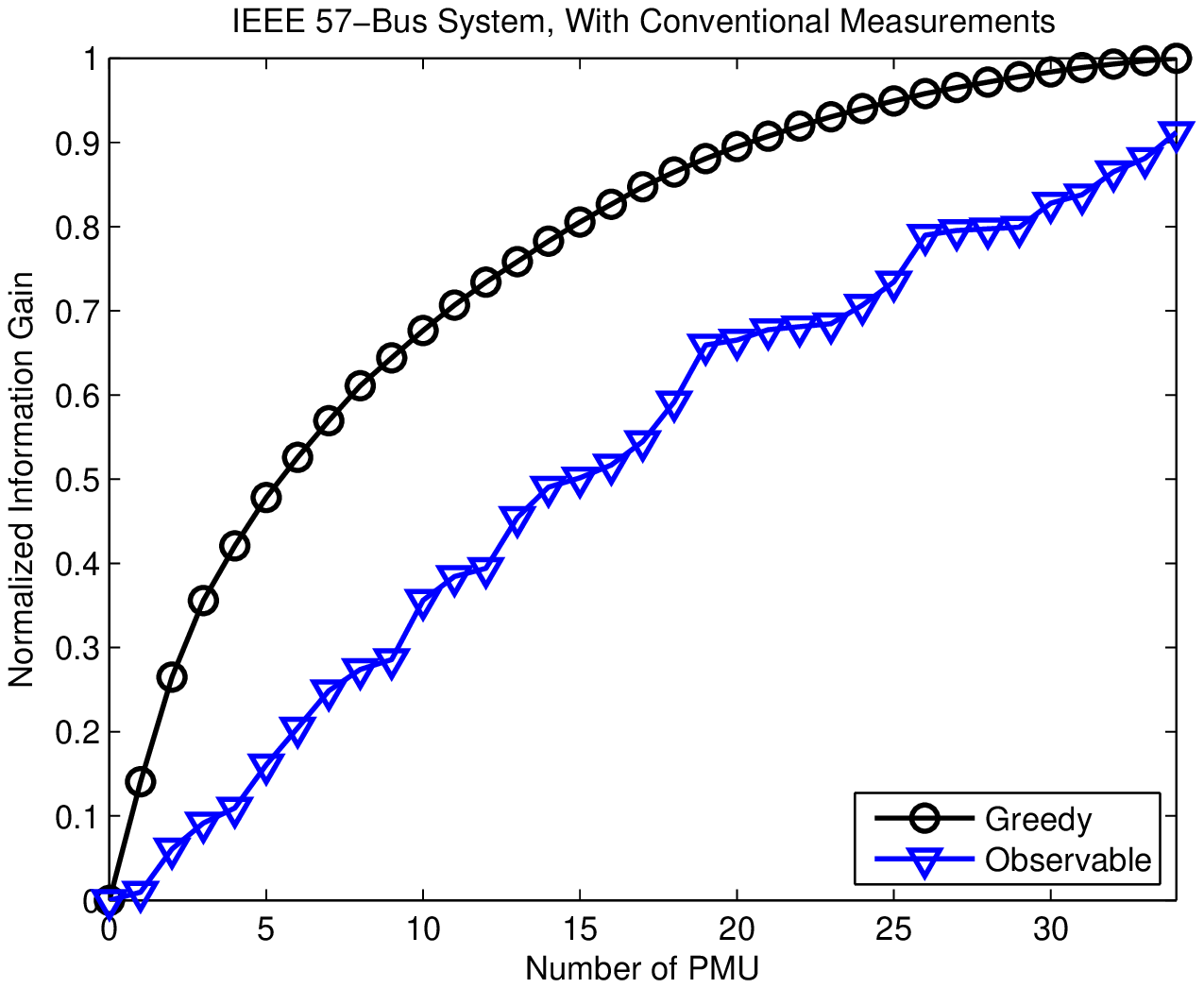,width=3.1in}\\
  (a) & (b)
  \end{tabular}
 \caption{Normalized information gain of different PMU placement schemes in the IEEE 57-bus system. (a) without conventional measurements and (b) with conventional measurements.}
 \label{fig:ieee118}
 \end{figure*}

\subsubsection{With Conventional Measurements}

In this case, real power flow measurements are assumed to be configured at all branches and buses. The resulting PMU configurations are shown in Table \ref{table:ieee14}. From the table, one can conclude that the optimal PMU placement is even more vulnerable to the changes in the PMU placement budget $K$ than the PMU only case, as the configurations changes significantly as $K$ increases. On the other hand, the greedy algorithm is robust, as $\ml S_g(K-1)\subset \ml S_g(K)$ for any $K\geq 1$. The normalized information gain is shown in Fig. \ref{fig:ieee14} (b). As the power system is already observable with conventional measurements, we use the same configuration for the `Observable' curve as the previous case. One can observe that the performance gain is larger compared to the case with only PMU measurements. This, again, confirms the conclusion that the pure topology based observability criterion can not efficiently model the uncertainties in the power system states.


%

\subsection{IEEE 57-Bus System}

\begin{table}[!t]
\renewcommand{\arraystretch}{1.3}
\caption{PMU Locations for 57-Bus System by Greedy Algorithm}
\label{table:ieee118}
\centering
\begin{tabular}{c|p{3.2cm}}
\hline
\bfseries Test Scenarios & \bfseries PMU Locations\\
\hline\hline
\begin{minipage}[t]{3.6cm}
  \centering
  IEEE 57-Bus (PMU Only)
\end{minipage}
& 9, 56, 18, 31, 12, 49, 29, 6, 25, 54, 20, 41, 38, 51, 32, 13, 27, 53, 57, 15, 19, 8, 30, 50, 17, 5, 16, 42, 52, 48, 55, 44, 24, 34 \\
\hline
\begin{minipage}[t]{3.6cm}
  \centering
 IEEE 57-Bus (Conventional)
\end{minipage}
& 56, 31, 19, 12, 54, 49, 25, 41, 32, 9, 29, 18, 6, 50, 20, 57, 27, 53, 38, 30, 13, 42, 51, 17, 55, 52, 5, 24, 34, 43, 16, 44, 8, 10\\
\hline
\end{tabular}
\end{table}

\subsubsection{PMU Measurements Only}

For the IEEE 57-bus system, it is computationally infeasible to obtain the optimal PMU configuration for large $K$. In such case, we only demonstrate the performance of the greedy PMU placement, and compare it against the observability based results in \cite{xu05}. For the case with only PMU measurements, the resulting PMU configurations are shown in Table \ref{table:ieee118}, and the normalized information gain is shown in Fig. \ref{fig:ieee118} (a). Similar to the IEEE 14-bus system, one can conclude that the greedy algorithm can achieve a significant information gain compared to the observability based criteria. This is because of the much more refined modeling of the MI function, which can effectively capture the remaining uncertainties in the states of the power system. Further, the information gain of the greedy placement curve also demonstrates the `diminishing marginal return' property.
 
\subsubsection{With Conventional Measurements}

In this case, real power flow measurements are assumed to be configured at all branches and buses. The resulting greedy PMU configurations are shown in Table \ref{table:ieee118}. The normalized information gain is shown in Fig. \ref{fig:ieee118} (b), where the greedy algorithm is compared against the same configuration in the previous case, Similar to the case without conventional measurements, one can observe the significant performance gain of the greedy algorithm. Note that, in this case, the power system is observable with conventional measurements. Therefore, the observability based methods fail to compute effective PMU configurations, due to the coarse modeling of the uncertainties in the power system states.

\section{Conclusion}
\label{sec:conclusion}

  This paper proposed an information-theoretic approach to address the phasor measurement units (PMUs)  placement for power system. Different from the topological observability based criterion in the literature, this paper proposed a much more refined, information-theoretic criterion, namely the mutual information (MI), as the PMU placement objective function. The proposed MI criterion can not only include the complete observability criterion as a special case, but also can accurately model the uncertainties in the system states. We further proposed a greedy PMU placement algorithm, and showed that it achieves an approximation ratio of $(1-1/e)$ for any PMU budget $K$, which is the best guarantee among polynomial-time algorithms. Such performance guarantee makes the greedy algorithm attractive in the typical scenario of phased installations, as the performance is robust to the changes in the PMU budget. Finally, the performance of the proposed PMU placement algorithm was demonstrated by simulation results.

\appendices
\section{Proof of Theorem \ref{thm:hardness}}
\label{apdx:hardness}
\begin{IEEEproof}
  We prove the hardness result by constructing polynomial-time reduction of an arbitrary instance of the max $k$-cover problem \cite{feige98} to a PMU placement problem. Thus, the hardness result easily follows from the $(1-1/e)$ inapproximability of the max $k$-cover problem \cite{feige98}. The max $k$-cover problem is as follows. We are given a set of elements $\ml U=\{1, 2, \ldots, N\}$, and a collection of sets $\ml P=\{\ml P_1, \ml P_2, \ldots, \ml P_M\}$, where each set $\ml P_i$ is a subset of $\ml U$. The task is to compute a subcollection of $K$ subsets $\ml P_{l_1}, \ml P_{l_2}, \ldots, \ml P_{l_K}$, such that the cardinality $|\cup_{k=1}^K \ml P_{l_k}|$ is maximized. Now, the reduction is as follows. Given a max $k$-cover problem instance, we construct a power system with $N+1$ buses, so that two buses $(i, j)$ are connected by a transmission line if and only if they both appear in a certain subset in $\ml P$. Further, we associate each subset $\ml P_i$ in $\ml P$ with a PMU, which is installed at a fixed (but can be arbitrarily chosen) bus in $\ml P_i$. Assume that the covariance matrix of the power injections $\bl P^{\text{inj}}$ is $\Sigma=\gamma B^2$, where $\gamma$ is a proper scalar that will be discussed later. Thus, according to Theorem \ref{thm:gmrf}, the resulting GMRF of $\bl \theta$ has covariance matrix $\gamma I$. Finally, choose $\gamma$ so that each discrete random variable $\theta_i$ has entropy 1. Assume there is no PMU failure, and zero PMU measurement noise, we can write the MI objective function as follows:
  \begin{eqnarray}
    I(\bl \theta; \bl z^{\text{PMU}}(\ml S))&=&H(\bl \theta)-H(\bl \theta|\bl z^{\text{PMU}}(\ml S))\\
    &\stackrel{(a)}{=}& N-H(\bl \theta|\bl z^{\text{PMU}}(\ml S))\\
    &\stackrel{(b)}{=}&N-|\overline{\cup_{k\in\ml S}\ml P_k}|\\
    &=&|{\cup_{k\in\ml S}\ml P_k}|
  \end{eqnarray}
  where (a) is because the phasor angles are independent, and each has entropy 1, by construction. (b) is because given the PMU measurements, the uncertainties only remain at the phasor angles of the `unobservable buses' $\overline{\cup_{k\in\ml S}\ml P_k}$. Denote $\ml P$ as the set of all possible PMU placement configurations. Thus, if we can solve the optimal PMU placement problem by maximizing $I(\bl \theta; \bl z^{\text{PMU}}(\ml S))$ subject to $|\ml S|\leq k$ beyond $(1-1/e)$ in polynomial time, we can also achieve the same performance guarantee for the max $k$-cover problem in polynomial time, from which the theorem follows.
\end{IEEEproof}

\section{Proof of Theorem \ref{thm:greedy}}
\label{apdx:greedy}

We now prove the $(1-1/e)$ performance guarantee of Algorithm \ref{alg:greedy}. The key lies in identifying the \emph{submodular} property of the PMU placement problem, which is a widely used concept in combinatorial optimizations.

\subsection{Introduction to Submodular Functions}

A set function $F$ is called submodular \cite{nemhauser78} if 
\begin{equation}
  F(\ml A\cup \{k\})-F(\ml A)\geq F(\ml B\cup \{k\})-F(\ml B)
  \label{eqn:submodular}
\end{equation}
for any sets $\ml A$ and $\ml B$ such that $\ml A\subseteq\ml B$. Essentially, this is the `diminishing marginal return' property, which, in the context of PMU placement, specifies that the marginal `information gain' is decreasing as the number of installed PMU increases. A set function $F$ is nondecreasing if $\ml S\subseteq \ml T$ implies that $F(\ml S)\leq F(\ml T)$, for all sets $\ml S$ and $\ml T$. The importance of submodularity can be seen by considering the following combinatorial optimization problem:
\begin{eqnarray}
  \ml S^\star=\arg\max_{|\ml S|\leq K}& F(\ml S)
  \label{eqn:opp_submodular}
\end{eqnarray}
For nondecreasing submodular functions, the following guarantee always holds for the greedy algorithm \cite{nemhauser78}:

\begin{lemma}
  Let a set function $F(\cdot)$ be submodular, nondecreasing and $F(\emptyset)=0$. For any $K\geq 1$, denote $\ml S^\star(K)$ and $\ml S_{g}(K)$ as the optimal solution to the problem (\ref{eqn:opp_submodular}) and the solution obtained by Algorithm \ref{alg:greedy}, respectively. Then, 
  \begin{equation}
    F(\ml S_{g}(K))\geq (1-{1/e})F(\ml S^\star(K))
  \end{equation}
  always holds. Thus, Algorithm \ref{alg:greedy} can achieve at least an approximation ratio of $(1-1/e)$.
  \label{lem:submodular_guarantee}
\end{lemma}

\subsection{Proof of Theorem \ref{thm:greedy}}

\begin{IEEEproof}
We now prove the theorem by showing that the objective functions $F_1(\cdot)$ and $F_2(\cdot)$ for the optimal PMU placement problem satisfy all assumptions in Lemma \ref{lem:submodular_guarantee}. We first fix time index $t$, and consider the case without conventional measurements. It is easy to see that $I_t(\bl \theta; \bl z^{\text{PMU}}(\emptyset))=0$. We next verify (\ref{eqn:submodular}) as follows:
\begin{eqnarray}
  &&I(\bl \theta; \bl z^{\text{PMU}}(\ml A\cup s))-I(\bl \theta; \bl z^{\text{PMU}}(\ml A))\\
  &\stackrel{(a)}{=}&I(\bl \theta; \bl z^{\text{PMU}}(s) |\bl z^{\text{PMU}}(\ml A))\label{eqn:mi_monotonic}\\
  &\stackrel{(b)}{=}&H(\bl z^{\text{PMU}}(s)|\bl z^{\text{PMU}}(\ml A))\nonumber\\
  &&\quad\qquad\qquad\qquad -H(\bl z^{\text{PMU}}(s)|\bl \theta, \bl z^{\text{PMU}}(\ml A))\\
  &\stackrel{(c)}{=}&H(\bl z^{\text{PMU}}(s)|\bl z^{\text{PMU}}(\ml A))-H(\bl z^{\text{PMU}}(s)|\bl \theta)
  \label{eqn:mi_submodular}
\end{eqnarray}
In above, time index is omitted for notation simplicity. $(a)$ follows from the chain rule of MI \cite{cover06}. $(b)$ follows from the definition of conditional MI. $(c)$ is because $\bl z^{\text{PMU}}(s)$ is conditionally independent of $\bl z^{\text{PMU}}(\ml A)$ given $\bl \theta$, due to the measurement model in (\ref{eqn:zpmui}) and (\ref{eqn:zpmuij}). This can also be clearly observed from the probabilistic graphical model in Fig. \ref{fig:sample}(b), where state variables $\bl \theta$ serve as the parent nodes of the PMU measurements in the `Bayesian part' of the graph. Note that $H(\bl z^{\text{PMU}}(s)|\bl z^{\text{PMU}}(\ml A))$ in (\ref{eqn:mi_submodular}) is decreasing in $\ml A$, as conditioning always reduces entropy \cite{cover06}. Since the second term in (\ref{eqn:mi_submodular}) is independent of $\ml A$, (\ref{eqn:submodular}) follows easily. Finally, as conditional MI is always nonnegative, one can deduce from (\ref{eqn:mi_monotonic}) that the MI function $I(\bl \theta; \bl z^{\text{PMU}}(\ml A))$ is also nondecreasing.

Now, the above analysis can be immediately extended to a time period of length $T$, where (\ref{eqn:submodular}) holds for $F_1(\cdot)$ by summing up the inequalities corresponding to each time slot $t$, and then dividing both sides by $T$. Thus, we conclude that the claim holds for $F_1(\cdot)$. Similarly, an identical analysis can be carried out in the case with conventional measurements, where all functions in (\ref{eqn:mi_monotonic})-(\ref{eqn:mi_submodular}) hold when conditioned on conventional measurements.
\end{IEEEproof}

\bibliographystyle{IEEEtran}
\bibliography{IEEEabrv,state_estimation}

%
%
%

%


%
%
%

    
\begin{IEEEbiographynophoto}{Qiao Li}
  (S'07) received the B.Engg. degree from the Department of Electronics Information Engineering, Tsinghua University, Beijing China, in 2006. He received the M.S. degree from the Department of Electrical and Computer Engineering, Carnegie Mellon University, Pittsburgh, PA USA, in 2008. He is currently a Ph.D. Candidate in the Department of Electrical and Computer Engineering, Carnegie Mellon University. His research interests include monitoring and control in electrical energy systems, cyber-physical systems, and wireless networking. 
\end{IEEEbiographynophoto}

\begin{IEEEbiographynophoto}{Tao Cui}
  (S'10) was born in Shaanxi Province, China, in 1985. He received the B.Sc. and M.Sc. degrees from Tsinghua University, Beijing, China, and is currently pursuing the Ph.D. degree at Carnegie Mellon University, Pittsburgh, PA. His main research interests include power system computation, protection, analysis, and control.
\end{IEEEbiographynophoto}

\begin{IEEEbiographynophoto}{Yang Weng}
  (S'08) received his B.Engg. in Electrical Engineering in 2006 from Huazhong University of Science and Technology, Wuhan, China. and a M.S. in Statistics from University of Illinois at Chicago in 2009. Currently, he is a graduate student in the Department of Electrical and Computer Engineering at Carnegie Mellon University, Pittsburgh, PA. His research interest includes power systems, smart grid, information theory, and wireless communications.
\end{IEEEbiographynophoto}

\begin{IEEEbiographynophoto}{Rohit Negi}
  (S'98-M'00) received the B.Tech. degree in electrical engineering from the Indian Institute of Technology, Bombay, in 1995. He received the M.S. and Ph.D. degrees from Stanford University, CA, in 1996 and 2000, respectively, both in electrical engineering. Since 2000, he has been with the Electrical and Computer Engineering Department, Carnegie Mellon University, Pittsburgh, PA, where he is a Professor. His research interests include signal processing, coding for communications systems, information theory, networking, cross-layer optimization, and sensor networks. Dr. Negi received the President of India Gold Medal in 1995.
\end{IEEEbiographynophoto}

\begin{IEEEbiographynophoto}{Franz Franchetti}
  received the Dipl.-Ing. degree and the PhD degree in technical mathematics from the Vienna University of Technology in 2000 and 2003, respectively. Dr. Franchetti has been with the Vienna University of Technology since 1997. He is currently an Assistant Research Professor with the Dept. of Electrical and Computer Engineering at Carnegie Mellon University. His research interests concentrate on the development of high performance DSP algorithms. 
\end{IEEEbiographynophoto}  

\begin{IEEEbiographynophoto}{Marija D. Ili\'c}
  (M'80-SM'86-F'99) is currently a Professor at Carnegie Mellon University, Pittsburgh, PA, with a joint appointment in the Electrical and Computer Engineering and Engineering and Public Policy Departments. She is also the Honorary Chaired Professor for Control of Future Electricity Network Operations at Delft University of Technology in Delft, The Netherlands. She was an assistant professor at Cornell University, Ithaca, NY, and tenured Associate Professor at the University of Illinois at Urbana-Champaign. She was then a Senior Research Scientist in Department of Electrical Engineering and Computer Science, Massachusetts Institute of Technology, Cambridge, from 1987 to 2002. She has 30 years of experience in teaching and research in the area of electrical power system modeling and control. Her main interest is in the systems aspects of operations, planning, and economics of the electric power industry. She has co-authored several books in her field of interest. Prof. Ili¡äc is an IEEE Distinguished Lecturer.
\end{IEEEbiographynophoto}

\end{document}